\documentclass[12pt]{article}
\usepackage{amsfonts,amssymb,amsmath}
\usepackage{amsfonts,amssymb}

\def\ZZ{{\mathbb Z}}

\def\Ok{{{\cal O}_K}}
\def\M{{M^*}}
\def\Pi{{\pi_*}}

\title{Notes on Abelian Class field theory}
\author{S. Subramanian}
\date{}

\begin{document}
\maketitle

\paragraph*{1.}  Let $K$ be a number field which for us would be a finite 
Galois extension of $Q$, the field of rational numbers (in particular, $Q$ 
itself is a number field).  The problem that is of  interest is to understand
Gal$(\overline{K}/K)^{ab}$ which is the abelianisation of 
Gal$(\overline{K}/K)$, the Galois group of $K$, where $\overline{K}$
denotes
an algebraic closure of $K$.  We let ${\cal O}_K$ denote the ring of 
integers of $K$, and $\ZZ$ the ring of integers of $Q$. 

Let $M$ be a finitely generated ${\cal O}_K$ submodule of $K$.  Since $M 
\subset K$, it is clear that $M \otimes_{{\cal O}_K} K = K$, so that rank 
of $M$ as an ${\cal O}_K$ module is one.  Let $M^*$ denote the dual 
${\cal O}$-module $M^* $ = Hom$_{{\cal O}_K}(M, {\cal O}_K)$.  Then 
$M^*$ is also a rank one $\Ok$-module, so that $M^*\otimes_{\Ok} K = K$. 
Since $M^*$ is finitely generated as an $\Ok$-module let $m_1, \cdots, m_r$
be generates for $M^*$.  The isomorphism $\M\otimes_{\Ok} K \to K$ enables 
us to regard $m_i \otimes 1$ as elements of $K$, so we see that $\M$ is also
a finitely generated $\Ok$ submodule of $K$.  It is now clear that 
$M \otimes_{\Ok}\M = \Ok$ and further $M \otimes_{\Ok}\M = M\M$ where on the 
right hand side, the multiplication is in $K$, regarding $M$ and $\M$ as 
submodules of $K$.  Let ${\cal C}\ell(K)$ denote the group of such finitely 
generated $\Ok$ submodules of $K$.  It is clear that ${\cal C}\ell(K)$ is an 
abelian group.

Let $p_1, p_2, q_1, q_2$ be prime elements of $\Ok$ (we emphasise : prime 
elements, not prime ideals) such that $p_1p_2 = q_1 q_2$ and the $p_i, q_i$ 
are all distinct.  Consider the $\Ok$-module $M$ generated by $\frac{1}{p_1},
\frac{1}{q_1}$, in $K$.  Then $M \otimes_{\Ok}~{\Ok}[\frac{1}{p_2}, 
\frac{1}{q_2}]$ 
is isomorphic to $\Ok [\frac{1}{p_2}, \frac{1}{q_2}]$ as an $\Ok[\frac{1}{p_2},
\frac{1}{q_2}]$ module, but $M$ is not isomorphic to $\Ok$ as an $\Ok$-module. 
This is a simple example of the fact that $\Ok$ is a UFD if and only if 
${\cal C}\ell(K) = 1$.  Since $K/Q$ is a finite Galois  extension, all but finitely
many primes in $\ZZ$ remain unramified in $\Ok$.  Let $\{p_1, \cdots,
p_m\}$ be
the set of primes in $\ZZ$ outside which Spec $\Ok \to  Spec \ZZ$ is 
unramified.  Let $S$ be the inverse image in $\Ok$ of the set $\{p_1, 
\cdots, p_m\}$.  We 
observe first that for $M \in {\cal C}\ell(K)$ we have an inclusion 
$$\Ok \subset M$$ 
such that $M/\Ok$ is a torsion $\Ok$-module .  Also, we have a strictly 
decreasing sequence 
$$M \supset M^2\supset M^3 \supset \cdots \supset \Ok$$
and hence it follows that $M^n = \Ok$ for some positive integer $n$, so that 
every element of ${\cal C}\ell(K)$ is of finite order. 

We need the following lemma:
\paragraph*{Lemma(1.1):} Let $X$ be an affine one-dimensional scheme (like 
Spec  $\ZZ$ or  Spec $\Ok$) and $\pi : Y \to X$ a finite Galois etale morphism.
Suppose every line bundle on $X$ is trivial.  Then every line bundle on $Y$ is
trivial.  
\paragraph*{Proof of Lemma (1.1):} Let $L$ be a line bundle on $Y$, possibly 
non trivial. We consider the vector bundle $\pi_* L$ on $X$.  Let rank $\pi_*L
= r$ = degree of $\pi$.  Since $X$ is affine and one dimensional, we obtain 
an exact sequence
$${\cal O} \to {\cal O}_X^{\oplus(r-1)} \to \pi_* L \to M \to 0$$ 
where $M$ is a line bundle on $X$ (equal to det $\Pi L$).  By hypothesis, 
$M$ is trivial, and on an affine scheme, extensions split, so $\Pi L$ is 
trivial.  This  implies that $L$ is trivial.  \hfill{Q.E.D.}

Let ${\cal O}_{K, S}$ denote the localisation of $\Ok$ obtained by
inverting all the elements of $S$, and $\ZZ_{\{p_1, \cdots, p_m\}}$ the
localisation of $\ZZ$ obtained by inverting $p_1, \cdots, p_m$.  Then the 
morphism $Spec {\cal O}_{K,S} \to Spec \ZZ_{\{p_1, \cdots p_m\}}$ is etale, 
and hence by Lemma 1 above every line bundle on Spec ${\cal O}_{K, S}$ is
trivial.  This forces:
\paragraph*{Lemma (1.2):}  Any $M \in {\cal C}\ell(K)$ satisfies $M
\subset
{\cal O}_{K, S}$.   

In particular:
\paragraph*{(Lemma 1.3:)} ${\cal C}\ell(K)$ is finite. 
\paragraph*{2.} Let $K$ be a number field as before, and $L/K$ a finite, 
Galois extension (not necessarily abelian), with Galois group $G$ and let 
$G^{ab}$ be the abelianisation of $G$, so that we have an exact sequence
$$1 \to [G, G] \to G \to G^{ab}\to 1.$$ 
We have 
\paragraph*{Proposition (2.1):} Let $M \in {\cal C}\ell(L)$ such that 
$M$ is not 
the pullback of an element of ${\cal C}\ell(K)$.  Then 
$$\sigma^*_1\sigma_2^* M =  \sigma^*_2 \sigma_1^* M$$ 
$\forall \sigma_1, \sigma_2 \in G$  such that $\sigma_1 \sigma_2 \not= 
\sigma_2 \sigma_1$.  
\paragraph*{Proof.}  Since $\sigma_1, \sigma_2$ do not commute in $G$, the 
orbit of $M$ under $<\sigma_1, \sigma_2>$ is a non-abelian subgroup of 
${\cal C}\ell(L)$ which is abelian (here  $<\sigma_1, \sigma_2)>$ 
denotes the group 
generated by $\sigma_1, \sigma_2)$.  It follows that the commutator $[G, G]$ 
acts trivially on ${\cal C}\ell(L)$. \hfill{Q.E.D.} 
\paragraph*{Proposition (2.2):} Any element $M \in {\cal C}\ell(L)$ 
fixed by 
$G^{ab}$, descends to an element of ${\cal C}\ell(K)$.  
\paragraph*{Proof:}  Follows from the above proposition and Galois descent.
\hfill{Q.E.D.} 
\paragraph*{Theorem (2.3):}  Let $K$ be a number field, and let $M\in 
{\cal C}\ell(K), M \not= \Ok$.  Let $M^n = {\cal O}$, where $n$ is the order of $M$.  
Then there is a finite cyclic $\ZZ/n$ extension $L/K$ such that $M$ becomes
trivial in ${\cal C}\ell(L)$.  
\paragraph*{Proof:} Let $\Ok$ be the ring of integers of $K$ and consider
the ring
$$R = \Ok \oplus M \oplus M^2 \oplus \cdots \oplus M^{n-1}$$ 
$R$ is an ${\cal O}_K$-algebra, and defines an integral extension of $\Ok$, 
whose  quotient field does the job. \hfill{Q.E.D.}
\paragraph*{Theorem (2.4):} Let $K$ be a number field such that ${\cal C}\ell(K)$ 
is nontrivial.  Then there exists a finite abelian Galois extension $L/K$ such
that every $M \in {\cal C}\ell(K)$ becomes the trivial element of ${\cal 
C}\ell(L)$. 
\paragraph*{Proof:}  By Theorem (2.3) above, we can do it for every element 
of ${\cal C}\ell(K)$, and since ${\cal C}\ell(K)$ is finite, we obtain a finite 
extension where this happens. \hfill{Q.E.D.}
\paragraph*{Theorem (2.5):}  Let $K$ be a number field such that ${\cal C}\ell(K)$ 
is nontrivial.  Then there is a finite, Galois, abelian extension $L/K$ whose 
Galois group is ${\cal C}\ell(K)$ such that every $M \in {\cal 
C}\ell(K)$ becomes 
trivial in ${\cal C}\ell(L)$.
\paragraph*{Proof:} Follows from previous steps. \hfill{Q.E.D.} 
\paragraph*{3.}  We now consider a number field $K$ such that 
${\cal C}\ell(K) = 1$.  As remarked before, it is easy to see that in this case
the ring of integers ${\cal O}_K$ is a UFD and hence a principal ideal domain. 
The typical case is $\ZZ$ in $Q$ and the arguments in the general case are 
similar. 

Let $L/Q$ be a finite, Galois, abelian extension of $Q$ with Galois group $G$.
Since $G$ is a finite abelian group, by the Chinese Remainder Theorem, 
$$G = \ZZ/_{{p_1}^{a_1}} \otimes \cdots \otimes \ZZ/_{{p_n}^{a_n}}$$ 
where $p_1, \cdots, p_n$ are rational primes.  By going modulo a subgroup 
of $G$ (every subgroup of $G$ is normal since $G$ is abelian) we may assume
that the Galois group of $L/K$ is $\ZZ/p^a$.  Unlike in the coprime case, when
we used the Chinese Remainder Theorem, and could have assumed the base field
was $Q$ without loss of generality, the group $\ZZ/p^a$ is a non split 
extension of $\ZZ/p$ factors.  We first consider the case when $L/K$ is a 
Galois extension with Galois group $\ZZ/p$, and as before, we consider the 
case $K= Q$ (the general case in similar).  Let ${\cal O}_L$ be the ring of 
integers of $L$ and let $q$ be a rational prime in $\ZZ$.  These are two cases 
to consider: $p\not= q, p = q$. 
\paragraph*{Case (i)} $q \not= p$. 
\paragraph*{Claim:} In this case, $Spec  {\cal O}_L \to Spec \ZZ$ is etale at 
$q$.  For, we consider a prime $q_1 \in {\cal O}_L$ such that $q^2_1$ divides 
$q$ in ${\cal O}_L$.  We consider the completion $L_{q_1}$ of $L$ at $q_1$ and 
the completion $Q_q$ of $Q$ at $q$.  We thus obtain an extension of local 
fields $L_{q_1}/Q_q$, again with Galois group $\ZZ/p$.  However, the residue
field extension ${\cal O}_L/q_1$ over $\ZZ/q$ is an extension of the finite
field $\ZZ/q$ and hence its Galois group is cyclic (generated by the Frobenius 
at $q$).  This cyclic group has to be a quotient of $\ZZ/p$ and hence has to 
be isomorphic to $\ZZ/p$.  This shows that $q$ remains unramified.  
\paragraph*{Case (ii) the case $q = p$.}   We recall that $L/Q$ is a Galois 
extension with Galois group $\ZZ/p$ and by Case (i) treated above, 
${\cal O}_L$  is unramified  outside $p$.  By Lemma (1.2) above, it follows 
that ${\cal C}\ell(L) =1$. 

From the above arguments, it follows that ${\pm}1 \in \ZZ$  
are the only points 
ramified in the extesion (possibly except for $p$) and hence the field 
extension is obtained by adjoining roots of unity.  

Further, since every time the class group remains trivial (in the case of 
a $\ZZ/p^a$ extension), we can repeat the argument.  We thus obtain 
\paragraph*{Theorem (3.1):}  Let $K$ be a number field with ${\cal C}\ell(K) = 1$. 
Then any abelian extension of $K$ is obtained by adjoining roots of unity. 
\paragraph*{Remark:}  The  exception occurs in the case of $\ZZ/2$ 
extension of $Q$, where $Q(\sqrt{p})$ is an abelian $\ZZ/2$ extension 
not obtained by adjoining a root of unity,where $p$ is a rational prime.  
This can be seen by looking at the arguments in Cases(i) and (ii) 
above.These fields have trivial class group by Lemma(1.2) above.  

\medskip
Address: School of Mathematics,Tata Institute of Fundamental 
Research,Mumbai400005,India\\

Email: subramnn@math.tifr.res.in

\end{document}